\documentclass[12pt]{amsart}
\usepackage{amssymb}
\def\C{\mathbb C}

\def\N{\mathbb N}

\begin{document}
\title[Functions of small growth]{Functions of small growth with no unbounded Fatou components}
\author{P.J. Rippon and G.M. Stallard}
\date{}
\maketitle

{\it Abstract.} We prove a form of the $\cos \pi \rho$ theorem
which gives strong estimates for the minimum modulus of a
transcendental entire function of order zero. We also prove a
generalisation of a result of Hinkkanen that gives a sufficient
condition for a transcendental entire function to have no
unbounded Fatou components. These two results enable us to show
that there is a large class of entire functions of order zero
which have no unbounded Fatou components. On the other hand we
give examples which show that there are in fact functions of
order zero which not only fail to satisfy Hinkkanen's condition
but also fail to satisfy our more general condition. We also give
a new regularity condition that is sufficient to ensure that a
transcendental entire function of order less than 1/2 has no
unbounded Fatou components. Finally, we observe that all the
conditions given here which guarantee that a transcendental
entire function has no unbounded Fatou components, also guarantee
that the escaping set is connected, thus answering a question of
Eremenko for such functions.

\section{Introduction}
\setcounter{equation}{0} Let $f$ be a transcendental entire
function and denote by $f^{n}, n \in \N$, the $n$th iterate of
$f$. The {\it Fatou set}, $F(f)$, is defined to be the set of
points, $z \in \C$, such that $(f^{n})_{n \in \N}$ forms a normal
family in some neighbourhood of $z$. The complement, $J(f)$, of
$F(f)$ is called the {\it Julia set} of $f$. An introduction to
the basic properties of these sets can be found in, for example,
[\ref{Berg}].\\

This paper concerns the question of whether a transcendental
entire function of small growth can have any unbounded components
of its Fatou set, a question which was first studied by Baker
[\ref{B81}]. An excellent survey article on this problem has been
written by Hinkkanen [\ref{H}].\\

The {\it order of growth}, $\rho$, of a transcendental entire
function $f$ is defined by
\[
  \rho = \overline{\lim_{r \to \infty}}\, \frac{\log \log M(r)}{\log
  r}, \mbox{ where } M(r) = \max_{|z| = r} |f(z)|.
\]
Defining $\sigma$ by
\[
 \sigma =  \overline{\lim_{r \to \infty}}\, \frac{\log
 M(r)}{r^{\rho}},
\]
we say that the growth of $f$ is of {\it minimal type} if $\sigma = 0$,
{\it mean type} if $0 < \sigma < \infty$ and {\it maximal type}
if $ \sigma = \infty$.\\

It was shown by Zheng [\ref{Z}, Theorem~1] that there are no
unbounded {\it periodic} Fatou components if the growth of $f$ is
at most order 1/2, minimal type. This result is sharp. For
example, Baker [\ref{B81}] showed that if $f(z) = z + \frac{\sin
\sqrt{z}}{\sqrt{z}} + c$, for $c$ sufficiently large, then $f$
has an unbounded invariant Fatou component (which is in
fact a Baker domain) and $f$ has order 1/2, mean type.\\

It is still unknown whether a transcendental entire function of order at most 1/2,
minimal type, can have any unbounded {\it wandering} components of
the Fatou set. Baker [\ref{B81}, Theorem~2] showed that there are
no unbounded Fatou components if
\begin{equation}
 \log M(r) = O((\log r)^p), \mbox{ for some } p \in (1,3).
\end{equation}
In [\ref{S}, Theorem~B] we strengthened this to show that there
are no unbounded Fatou components if
 there exist $\epsilon \in (0,1)$ and $R>0$ such that
\begin{equation}
   \log \log M(r) < \frac{(\log r)^{1/2}}{(\log \log
   r)^{\epsilon}}, \; \mbox{ for } r>R.
\end{equation}

Although many authors have worked on this problem, all further
results of this type have required some regularity condition on the growth. We
discuss some of these regularity conditions in Section~7.\\

Hinkkanen [\ref{H2}, Theorem~1] showed that a transcendental
entire function of order less than 1/2 has no unbounded Fatou
components if there exist $0 < \delta \leq 1$, $L > 1$ and $C>0$
such that, for large $r$, there exists $t \in (r, r^L)$ with
\begin{equation}
 \frac{\log m(t)}{\log M(r)} \geq L \left( 1 - \frac{C}{(\log
 r)^\delta} \right),
\end{equation}
where $m(r) = \min _{|z| = r} |f(z)|$.\\

He also suggested that it was plausible that condition (1.3) is
satisfied by all functions of order at most 1/2, minimal type. We shall show that there are in fact functions of order zero which fail to satisfy (1.3).\\

First, however, in Section~2 we obtain the following generalisation of
Hinkkanen's result, using a somewhat simpler method of proof. In
Section~5 we show that the conditions of this new result are
satisfied for a much larger class of functions of order zero than
those which satisfy (1.2). It follows that such functions have no
unbounded Fatou components.

\newtheorem{1}{Theorem}
\begin{1}
Let $f$ be a transcendental entire function and let $R_1 > 0$ be
such that the sequence $R_n$ defined by $R_{n+1} = M(R_n)$ tends to $\infty$. Suppose that there exist $L>1$ and a positive
sequence $a_n$ such that, for all $r \in [R_n,R_{n+1})$, there
exists $t \in (r,r^L)$ with
\begin{equation}
\frac{\log m(t)}{\log M(r)} \geq L(1 - a_n)
\end{equation}
and
\[
\sum_{n \in \N} a_n < \infty.
\]
Then $F(f)$ has no unbounded components.
\end{1}

{\it Remark.} Hinkkanen's result follows from Theorem~1 by taking
$a_n = C/(\log R_n)^{\delta}$ for some $C>0$, $0 < \delta \leq
1$, where $R_{n+1} = M(R_n)$. To see that this is true, note that,
for any transcendental entire function $f$, if $R_n$ tends to
$\infty$ then
\[
\log M(R_n) / \log R_n \to \infty\; \mbox{ as } n
\to \infty.
\]
(In what follows we use this fact
without further reference.) Thus
\begin{equation}
  \log R_n > 2^n, \mbox{ for large }n,
\end{equation}
and so $\sum_{n \in \N} a_n  = \sum_{n \in \N} \frac{C}{(\log
R_n)^{\delta}} < \infty$. We observe that Theorem~1
does not require the assumption that the order of $f$ is less than $1/2$.\\

Many of the results on this subject use the version of the $\cos
\pi \rho$ theorem proved by Barry [\ref{Ba}, p.294]. This states
that if $f$ has order $\rho < 1/2$ then, for each $\alpha \in
(\rho, 1/2)$, the set of values of $r$ for which
\[ \log m(r) > \cos \pi \alpha \log M(r)
\]
has lower logarithmic density at least $1 - \rho/\alpha$. In order
to show that the hypotheses of Theorem~1 are satisfied for many
functions of order zero, in Sections~3 and 4 we prove the following result which can be regarded as a
version of the
$\cos \pi \rho$ theorem.

\newtheorem{6}[1]{Theorem}
\begin{6}
Let $f$ be a transcendental entire function of order $\rho$, $\rho
<\alpha < 1$ and put
\[
  \epsilon(r) = \frac{\log \log M(r)}{\log r} \; \mbox{ and } \;
  \delta(r) = \max_{r \leq t \leq r^{1/(1-\alpha)}} \epsilon(t).
\]
Let $0 <\eta < 1/2$ and $\mu,$ $\nu$ be functions such that,
for large $r$,
\[
0 < \mu(r) \leq 1,\; 0 < \nu(r) \leq 1/4 \; \mbox{ and } \; \mu(r) \nu(r) >
\frac{2 \epsilon(r)}{1 - \delta(r)}.
\]
Then, for large $r$, there exists $R \in (r^{(1 - \delta(r))(1 -
\mu(r))},r^{1 - \delta(r)})$ such that
\[
 \log m(t) > (1 - 20 \log (2e/\eta) \nu(r)) \log M(t),
\]
for $0 \leq t \leq R/2$ except in a set of intervals, the sum of
whose lengths is at most $\eta R$.
\end{6}

Theorem~2 gives a big improvement on existing estimates for the
minimum modulus near values of $r$ for which $M(r)$ is small. In
particular, it gives very precise information for functions of
order zero. In Section~5 we apply this result with $\rho = 0$,
$\mu(r) = \sqrt{\epsilon(r)}$ and $\nu(r) = 3 \sqrt{\epsilon(r)}$,
in which case $\mu(r) \to 0$ and $\nu(r) \to 0$ as $r \to \infty$.
Together with Theorem~1, this enables us to prove the following
result.

\newtheorem{8}[1]{Theorem}
\begin{8}
Let $f$ be a transcendental entire function of order zero and let
$R_1>0$ be such that the sequence $R_n$ defined by $R_{n+1} =
M(R_n)$ tends to $\infty$ and the sequence $\epsilon_n$
defined by
\[
  \epsilon_n = \max_{R_n \leq r \leq R_{n+1}} \frac{\log \log M(r)}{\log
  r}
\]
is positive. If
\[
  \sum_{n \in \N} \sqrt{\epsilon_n} < \infty,
\]
then $F(f)$ has no unbounded components.
\end{8}

We conclude Section~5 by using Theorem~3 to prove the following
result which relates the size of $M(r)$ to the existence of unbounded Fatou components.\\

{\it Remark.} By using a more sophisticated proof, we can remove the square root in the
condition in Theorem~3. This improvement in Theorem~3 does not,
however, lead to an improvement in Theorem~4 below and so we omit the details.

\newtheorem{9}[1]{Theorem}
\begin{9}
Let $f$ be a transcendental entire function and suppose that
 there exist $m \in \N$ and $R>0$ such that
\begin{equation}
   \log \log M(r) < \frac{\log r}{ \log^m r},
    \; \mbox{ for } r>R.
\end{equation}
 Then $F(f)$ has no unbounded components.
\end{9}

Note that $\log ^m$ denotes the $m$th iterated logarithm function.\\

Condition (1.6) is a significant improvement on condition (1.2)
which was previously the best condition of this type. It still
does not, however, cover all functions of order zero which need
only satisfy
\[
  \log \log M(r) = o(1)\log r \; \mbox{ as } r \to \infty.
\]

In Section~6, we study functions of the form
\begin{equation}
  f(z) = \prod_{m=1}^{\infty} \left(1 -
  \frac{z}{r_m}\right)^{r_m^{\epsilon_m}},
\end{equation}
where $\epsilon_m$ is a decreasing null sequence  and $r_m$ is a
strictly increasing sequence such that $r_m^{\epsilon_m} \in \N$ for $m \in \N$.
We show that if $\epsilon_m$ and
$r_m$ satisfy certain conditions, then $f$ is a function of order
zero that fails to satisfy Hinkkanen's condition (1.3). More
generally, we show that there are functions of the form (1.7)
that have order zero and fail to satisfy the hypotheses of
Theorem~1. These examples suggest that new techniques are needed
to answer Baker's original question.\\

As mentioned earlier, many authors have shown that a transcendental entire function of
order less than 1/2 has no unbounded Fatou components provided
that the growth of the function is sufficiently regular. In
Section~7 we prove the following result of this type.

\newtheorem{10}[1]{Theorem}
\begin{10}
Let $f$ be a transcendental entire function of order $\rho <
1/2$. Suppose that, for each $m>1$, there exists a real function
$\psi$ defined on $(r_0, \infty)$, where $r_0>0$, such that, for
$r\geq r_0$, we have $\psi(r) \geq r$ and
\begin{equation}
  M(\psi(r)) \geq \psi(M(r))^m.
\end{equation}
Then $F(f)$ has no unbounded components.
\end{10}

We show that many of the known regularity conditions that have
been used in connection with Baker's question can be written in
the form (1.8) for some function $\psi$. We end Section~7 by
proving the following result which gives a new regularity
condition that is sufficient for the hypotheses of Theorem~5 to
be satisfied. Here $\exp^n$ and $\log^n$ denote the $n$th
iterated exponential and logarithm functions, respectively.

\newtheorem{11}[1]{Theorem}
\begin{11}
Let $f$ be a transcendental entire function of order $\rho < 1/2$
and suppose that there exist $n \in \N$ and $0<q<1$  such that
\begin{equation}
 M(r) \geq \exp^{n+1}((\log^nr)^q), \; \mbox{ for large } r.
\end{equation}
Then the hypotheses of Theorem~5 hold with
\[
 \psi(r) = \exp^{n}((\log r)^p), \; \mbox{ where } pq>1.
\]
Hence $F(f)$ has no unbounded components.
\end{11}

We remark that condition (1.9) becomes less restrictive as $n$
increases. Indeed, if $0<q<1$ and $n \in \N$, then there exists
$r(n,q)>0$ such that
\[
  \exp^{n+2}((\log^{n+1}r)^q) \leq \exp^{n+1}((\log^nr)^q), \mbox{
  for } r \geq r(n,q).
\]
There are, however, many transcendental entire functions of order
$\rho < 1/2$ that fail to satisfy (1.9) for any $n \in \N$ and
$0<q<1$. Indeed, (1.9) implies that $\phi(t) = \log M(e^t)$
satisfies
$\phi(t)/t^k \to \infty$ as $t \to \infty$, for all $k>1$.\\

Theorem~6 generalises a result of Singh and Taniguchi (see, for
example, [\ref{Si}]). They showed that a transcendental entire
function of order $\rho < 1/2$ satisfying (1.9) with $n=1$ has no
unbounded Fatou components. Their result in turn generalised a
result of Wang [\ref{W}] who considered functions satisfying
(1.9) with $n=0$, that is, functions of positive lower order; an alternative proof of Wang's result
is given in [\ref{H}, p.205].\\

Finally, we point out that our results give several new sufficient conditions for the escaping set of a transcendental entire function to be connected. In [\ref{RS}, Corollary 5] we show that if $f$ is a transcendental entire function
 that satisfies (1.2) or has order $\rho < 1/2$ and satisfies a certain regularity condition, then the escaping set,
 \[
  I(f) = \{z: f^n(z) \to \infty\; \mbox{ as } n \to \infty\},
 \]
 is connected.  We remarked in that paper that in fact the escaping set is an unbounded connected set for any
 function satisfying the hypotheses of Lemma~2.1 below. Thus such functions satisfy Eremenko's conjecture [\ref{E}]
 that the escaping set of a transcendental entire function has no bounded components. In this paper, we show
 that the hypotheses of Lemma~2.1 are satisfied by any function satisfying the hypotheses of Theorem~1, Theorem~3,
 Theorem~4, Theorem~5 or Theorem~6. Therefore, in each of these theorems, the conclusion that $F(f)$ has no unbounded
 components can be replaced by the conclusion that $I(f)$ is connected and hence Eremenko's conjecture holds in these cases.\\

 Note that in [\ref{RS}], when showing that the hypotheses of Lemma~2.1 are sufficient to ensure that the escaping set is connected, we actually obtain more detailed information about the structure of the escaping set for such functions. More precisely, we consider the set of fast escaping points defined by
\[
  B(f) = \{z: \mbox{ there exists } L \in \N \mbox{ such that
  } f^{n+L}(z) \notin \widetilde{f^n(D)}, \mbox{ for } n \in \N \}
\]
and the subset of $B(f)$ defined by
 \[
  B_D(f) = \{ z: f^n(z)\notin \widetilde{f^n(D)}, \mbox{ for } n \in \N \},
\]
where $D$ is any open disc meeting $J(f)$ and $\widetilde{U}$
denotes the union of $U$ and its bounded complementary
components. (Note that the set $B(f)$ is independent of the choice
of $D$. The set $B(f)$ was first introduced in [\ref{RS05}] where
we showed that it is equal to the set $A(f)$ introduced by
Bergweiler and Hinkkanen in [\ref{BH}].) We show in [\ref{RS}]
that if a transcendental entire function $f$ satisfies the
hypotheses of Lemma~2.1,
then $B_D(f)^c$ has a bounded component from which it follows that $B_D(f)$, $B(f)$ and $I(f)$ are all connected.\\

{\it Remark.} After preparing this paper, we learnt that results
closely related to Theorem~1, Theorem~3, Theorem~4 and the
examples in Section~6 have also been obtained by Hinkkanen and
Miles; see [\ref{HM}].

\section{Proof of Theorem~1}
\setcounter{equation}{0}

Many of the papers on this subject use the following result
(proved in [\ref{S}, Lemma~2.7]) which is a generalisation of a
result by Baker [\ref{B81}, proof of Theorem~2].

\newtheorem{10.1}{Lemma}[section]
\begin{10.1}
Let $f$ be a transcendental entire function and suppose that
there exist sequences $R_n, \rho_n \to \infty$ and $c(n) > 1$
such that
\begin{enumerate}
\item $R_{n+1} = M(R_n)$,
\item $R_n \leq \rho_n \leq R_n^{c(n)}$,
\item $m(\rho_n) \geq R_{n+1}^{c(n+1)}$.
\end{enumerate}
Then $F(f)$ has no unbounded components.
\end{10.1}

Our proof of Theorem~1 is based on Lemma~2.1. We also use the
following consequence of the convexity of $\log M(e^t)$.

\newtheorem{10.2}[10.1]{Lemma}
\begin{10.2}
Let $f$ be a transcendental entire function. Then there exists
$R>0$ such that, for all $r \geq R$ and all $c>1$,
\[
  \log M(r^c) \geq c  \log M(r).
\]
\end{10.2}

\begin{proof} Since $\phi(t) = \log M(e^t)$ is convex and
$\phi(t)/t \to \infty$ as $t \to \infty$, we have
\[
  \phi'(t) \geq \frac{\phi(t)}{t}, \mbox{ for } t \geq T,
\]
say. Here, for definiteness, $\phi'$ denotes the right derivative of $\phi$. Thus, for $t \geq T$ and any $c>1$,
\[
 \log \frac{\phi(ct)}{\phi(t)} = \int_t^{ct}
 \frac{\phi'(u)}{\phi(u)}du
 \geq \int_t^{ct} \frac{1}{u}du = \log c.
\]
Hence the result holds with $R=e^T$.
\end{proof}

{\it Proof of Theorem~1.}\;
Let $f$ be a function satisfying the
hypotheses of Theorem~1 and suppose that $\sum_{n \in \N} a_n <
\infty$.  We may assume that
 $a_n < 1$, for all $n \in \N$. Since $\sum_{n \in \N}  a_n  < \infty$, we have
\[
  0 < \prod_{n \in \N} \left( 1 - a_n  \right) <
  1
\]
and so, for a given $L>1$, we can define
\[
 c(n) = \frac{L}{\prod_{m=n}^{\infty} (1 - a_m)}, \mbox{ for } n \in \N.
\]
Note that
\begin{equation}
 L < c(n) < \infty \; \mbox{ and } \; c(n+1) = \left( 1 - a_n \right) c(n) < c(n).
\end{equation}
We will show that the conditions of Lemma~2.1 are satisfied for
these values of $c(n)$. We may assume that $R_1$ has been chosen
sufficiently large to ensure that $R_n^{c(n)/L} < R_{n+1}$ for
each $n \in \N$. So, it follows from (1.4), (2.1) and Lemma~2.2
that, for each $n \in \N$, there exists $\rho_n \in (R_n^{c(n)/L},
R_n^{c(n)}) \subset (R_n, R_n^{c(n)})$ such that
\begin{eqnarray*}
\log m(\rho_n) & \geq & L(1 - a_n) \log
M(R_n^{c(n)/L})\\
& \geq & L(1 - a_n) \frac{c(n)}{L}
\log M(R_n)\\
& = & c(n+1) \log M(R_n) = c(n+1) \log R_{n+1}.
\end{eqnarray*}
Thus, for all $n \in \N$, there exists $\rho_n \in (R_n,
R_n^{c(n)})$ such that $m(\rho_n) \geq R_{n+1}^{c(n+1)}$. The
result of Theorem~1 now follows from Lemma~2.1.

\section{Results of Cartwright}
\setcounter{equation}{0}

This section and the next are concerned with the proof of
Theorem~2. We first consider the case that $f(0) = 1$
and use the following notation:\\

$n(r)$ is the number of zeros of $f$ in $\{z: |z| \leq r\}$, counted according to multiplicity,\\

$\displaystyle N(r) = \int_0^r \frac{n(t)}{t}\,dt$,\\

$\displaystyle Q(r) = r \int_r^{\infty} \frac{n(t)}{t^2}\, dt$,\\

$\displaystyle B(r) = r \int_0^{\infty} \frac{n(t)}{t(r+t)}\,dt$,\\

$\displaystyle a(r) = r B'(r) = r \int_0^{\infty} \frac{n(t)}{(r+t)^2}\,dt$.\\

It is well known that these integrals are convergent if $f$ has order $\rho < 1$; this follows from Lemma~3.1(b) below.
We show later that $B(r) = \log M(r)$ in the case that all the
zeros of $f$ are on the negative real axis.\\

 Cartwright [\ref{C},
Theorem~51 and Theorem~52] proved careful estimates for the
minimum and maximum modulus of a function of order zero in terms
of the quantities $N(r)$ and $Q(r)$. She then used these
estimates to show that, for a function of order zero, $\log m(r)$
is asymptotically equal to $\log M(r)$ in a set of upper density
1. In this paper we build upon her results to show how the ratio
$\log m(r)/\log M(r)$ depends on the rate of growth of the
function.\\

 The key idea in the proof is to use the quantities
$a(r)$ and $B(r)$ to estimate
$Q(r)$ and $N(r)$. An advantage of
considering $a(r)/B(r)$ is that $a(r)$ is defined in terms of the
derivative of $B(r)$. (Alternatively, similar arguments can be used with the quantities $Q(r)$ and $N(r) + Q(r)$
since $Q(r)$ can be expressed in terms of the derivative of $N(r) + Q(r)$. Each approach has some
advantages and some disadvantages.)\\

In this section we give some preliminary results that will be
useful in the proof of Theorem~2. We begin by noting various
relationships between the quantities defined above.

\newtheorem{2.1}{Lemma}[section]
\begin{2.1}
Let $f$ be a transcendental entire function of order $\rho < 1$ with $f(0) = 1$.
Then, for $r>0$,\\

(a) $N(r) < \log M(r)$,\\

(b) $n(r) < \log M(3r)$,\\

(c) $n(r) \leq Q(r)$,\\

(d) $Q(r) \leq 4a(r)$,\\

(e) $B(r) \leq N(r) + 2a(r)$ and $B(r) \leq N(r) + Q(r)$.
\end{2.1}
\begin{proof}
{\it (a)} This follows from Jensen's theorem (see, for example,
[\ref{T}, second
edition, p.125] for a proof).\\

{\it (b)} It follows from (a) that
\[
 n(r) \log 3 \leq \int_r^{3r} \frac{n(t)}{t}dt \leq N(3r) < \log
 M(3r).
\]

{\it (c)} We have
\[
Q(r) =  r \int_r^{\infty} \frac{n(t)}{t^2} dt \geq r n(r)
\int_r^{\infty} \frac{1}{t^2} dt = n(r).
\]

{\it (d)} We have
\[
Q(r) = r \int_r^{\infty} \frac{n(t)}{t^2} dt
  \leq 4r \int_r^{\infty} \frac{n(t)}{(r+t)^2} dt
  \leq 4a(r).
\]

{\it (e)} We have
\[
B(r)  =  r
\int_0^{\infty}\frac{n(t)}{t(r+t)}dt
  =  r \int_0^r\frac{n(t)}{t(r+t)}dt + r
 \int_r^{\infty}\frac{n(t)}{t(r+t)}dt.
 \]
 Thus
 \[
 B(r) \leq  N(r) + 2a(r)\; \mbox{ and also }\; B(r) \leq N(r) + Q(r),
 \]
 as required.
 \end{proof}

 The next result shows that, for a function of order less than 1,
  the growth of $B(r)$ is limited by the
 growth of $M(r)$.

\newtheorem{2.2}[2.1]{Lemma}
\begin{2.2}
Let $f$ be a transcendental entire function of order $\rho$, $\rho
<\alpha < 1$ with $f(0) = 1$ and put
\[
  \epsilon(r) = \frac{\log \log M(r)}{\log r} \; \mbox{ and } \;
  \delta(r) = \max_{r \leq t \leq r^{1/(1-\alpha)}} \epsilon(t).
\]
Then, for large $r$, we have $0 < \delta(r) < \alpha < 1$ and
\[
 B(r^{1 - \delta(r)}) < 3 \log M(r) \log r.
\]
\end{2.2}

\begin{proof} Clearly $0 < \delta(r) < \alpha < 1$ for large $r$, since $\rho < \alpha$. It follows from Lemma~3.1(b) that, for large $t$,
\[
 n(t) < \log M(3t) < (3t)^{\alpha}.
\]
Also, for large $r$ and $r/3 \leq t \leq r^{1/(1-\alpha)}/3$,
\[
 n(t) < \log M(3t) < (3t)^{\delta(r)}.
\]

 Using these estimates for
$n(t)$ and the fact that $f(0) = 1$ so that $n(t) = 0$ for $0 \leq
t \leq r_1$, say, we have, for large $r$,
\begin{eqnarray*}
B(r^{1 - \delta(r)}) & = & r^{1 - \delta(r)} \int_0^{r/3}
\frac{n(t)}{t(r^{1 - \delta(r)}+t)}dt + r^{1 - \delta(r)}
\int_{r/3}^{r^{1/(1 - \alpha)}/3} \frac{n(t)}{t(r^{1 -
\delta(r)}+t)}dt\\
&&+ \, r^{1 - \delta(r)} \int_{r^{1/(1 - \alpha)}/3}^{\infty} \frac{n(t)}{t(r^{1 - \delta(r)}+t)}dt \\
& \leq & \log M(r) \int_{r_1}^{r/3} \frac{1}{t} dt + r^{1
-\delta(r)}\int_{r/3}^{r^{1/(1 - \alpha)}/3}
\frac{(3t)^{\delta(r)}}{t^2}dt\\
&& + \, r^{1 - \delta(r)}\int_{r^{1/(1 -
\alpha)}/3}^{\infty} \frac{(3t)^{\alpha}}{t^2}dt\\
& < & 2\log M(r) \log (r/3) + \frac{3^{\delta(r)}r^{1 - \delta(r)}}{(r/3)^{1 - \delta(r)}(1 - \delta(r))}
 + \frac{3^{\alpha} r^{1 - \delta(r)}}{(r^{1/(1-\alpha)}/3)^{1-\alpha}(1-\alpha)}\\
& < & 2\log M(r) \log (r/3) + \frac{3}{1 - \delta(r)} +
\frac{3}{1 - \alpha}\\
& < & 3 \log M(r) \log r,
\end{eqnarray*}
since $\delta(r) < \alpha$ for large $r$.
\end{proof}

The next result shows that, for a function of order less than 1,
$\log M(r)$ is always bounded above by $B(r)$. This is similar to
Cartwright's result [\ref{C}, Theorem~51]. Note that the proof
shows that we have equality in the case that all the zeros of $f$
are on the negative real axis. This result follows, for example,
from arguments in [\ref{T}, p.271] but we include it here for the
sake of completeness.

\newtheorem{2.3}[2.1]{Lemma}
\begin{2.3}
Let $f$ be a transcendental entire function of order less than $1$
with $f(0) = 1$. Then, for $r > 0$,
\[
 \log M(r) \leq B(r) \leq N(r) + Q(r).
\]
\end{2.3}

\begin{proof} Let $f$ be a transcendental entire function of order
less than 1 with $f(0) = 1$ and let the zeros of $f$ be at the
points $z_n$, $n \in \N$, with modulus $r_n$, where $0 < r_1 \leq r_2 \leq \cdots$. It follows from
Hadamard's factorization theorem (see, for example, [\ref{T},
p.250]) that
\[
  f(z) = \prod_{n \in \N} (1 - z/z_n)
\]
and so
\[
 |f(z)| \leq \prod_{n \in \N} (1 + |z|/r_n).
\]
Thus
\begin{eqnarray*}
 \log M(r) & \leq & \sum_{n \in \N} \log (1 + r/r_n)
 = \sum_{n \in \N} n \left( \log (1 + r/r_n) - \log (1 +
 r/r_{n+1})
 \right)\\
 & = & \sum_{n \in \N} n \int_{r_n}^{r_{n+1}} \frac{r}{t(r+t)}dt =
 B(r).
\end{eqnarray*}
The result now follows from Lemma~3.1(e).
\end{proof}

It is much harder to obtain estimates for the minimum modulus of
an entire function. We will use the careful estimates obtained by
Cartwright in [\ref{C}, Theorem~52]. Her results are stated only
for functions of order zero but the proofs apply for functions of
order less than 1. Again, since this result is crucial to our
argument, we include a proof for the sake of completeness. We
state our results just for the case that $f(0) = 1$. The proof
follows that of Cartwright which depends on the following
covering lemma of Cartan [\ref{Ca}, p.273].

\newtheorem{2.4}[2.1]{Lemma}
\begin{2.4}
Let $z_1, \ldots, z_m$ be $m$ points distinct or coincident in
$\C$, and let $h>0$. Then $\{z: \prod_{n=1}^m |z-z_n| \leq h^m\}$
can be enclosed in at most $m$ discs, the sum of whose radii is
at most $2eh$.
\end{2.4}

Following Cartwright, we now obtain a lower estimate for the
minimum modulus in terms of $N(r)$ and $Q(r)$.

\newtheorem{2.5}[2.1]{Lemma}
\begin{2.5}
Let $f$ be a transcendental entire function of order less than $1$
with $f(0) = 1$, and let $0 < \eta < 1/2$. Then, for large $R$,
\[
  \log m(r) > N(R) - (1 + \log(2e/\eta)) Q(R),
\]
for $0 \leq r \leq R/2$ except in a set of intervals, the sum of
whose lengths is at most $\eta R$.
\end{2.5}

\begin{proof} As in the proof of Lemma~3.3, we write
\[
  f(z) = \prod_{n \in \N} (1 - z/z_n),
\]
where the points $z_n$ are the zeros of $f$ and have modulus
$r_n$, with $0 < r_1 \leq r_2 \leq \cdots$. Suppose that $|z| = r
\leq R/2$, where $R>r_1$, and let $m$ denote the largest integer
for which $r_m < R$. Then
\begin{equation}
 \log |f(z)| = \log \prod_{n =1}^m |1 - z/z_n| + \log \prod_{n = m+1}^{\infty} |1 -
 z/z_n|.
\end{equation}

We first obtain a lower bound for the last term in (3.1). Note
that
\begin{equation}
  \log \prod_{n = m+1}^{\infty} |1 -
 z/z_n| \geq \log \prod_{n = m+1}^{\infty} (1 -
 r/r_n) = \sum_{n = m+1}^{\infty} \log (1- r/r_n).
\end{equation}

If $n>m$ then $r_n \geq R \geq 2r$ and so $0 < r/r_n \leq 1/2$. Thus
\[
\log (1- r/r_n) > -2r/r_n
\]
and so
\begin{equation}
 \sum_{n= m+1}^{\infty} \log (1 - r/r_n) > - \sum_{n=m+1}^{\infty}
 2r/r_n.
\end{equation}
Now
\begin{eqnarray*}
  \sum_{n= m+1}^{\infty} r/r_n & \leq & r\sum_{n= m+1}^{\infty}n( 1/r_n -
  1/r_{n+1})\\
  &  = & r \sum_{n= m+1}^{\infty} n\int_{r_n}^{r_{n+1}}
  \frac{1}{t^2}dt \\
  & \leq & r \int_{R}^\infty \frac{n(t)}{t^2}dt \leq
  \frac{1}{2}Q(R).
\end{eqnarray*}
Thus, by (3.2) and (3.3),
\begin{equation}
 \log \prod_{n = m+1}^{\infty} |1 -
 z/z_n| \geq \sum_{n= m+1}^{\infty} \log (1 - r/r_n) \geq - Q(R).
\end{equation}

We now estimate $\log \prod_{n =1}^m |1 - z/z_n|$. To do this we
write $z = Rz'$ (so that $|z'| \leq 1/2$) and note that
\[
|1 - z/z_n| = |1 - Rz'/z_n| = \left| \frac{R}{z_n}\right| \left|
\frac{z_n}{R} - z'\right|.
\]
Thus
\begin{equation}
  \log \prod_{n =1}^m |1 - z/z_n| =
\log \prod_{n =1}^m \frac{R}{r_n} + \log \prod_{n =1}^m \left|
\frac{z_n}{R} - z'\right|.
\end{equation}

We estimate the last term in (3.5) using Cartan's lemma (Lemma
3.4) with $h = \eta/(2e)$. This implies that
\begin{equation}
   \prod_{n =1}^m
\left| \frac{z_n}{R} - z'\right| > \left(\frac{\eta}{2e}\right)^m,
\end{equation}
except for $z'$ in a set of discs, the sum of whose radii is at
most $\eta$; that is, except for $z$ in a set of discs, the sum
of whose radii is at most $\eta R$. Thus, by Lemma~3.1(c),

\[
  \log \prod_{n =1}^m
\left| \frac{z_n}{R} - z'\right| > m \log (\eta/(2e)) = -m \log
(2e/\eta)
\]
\begin{equation}
\geq -n(R) \log (2e/\eta) \geq -Q(R) \log (2e/\eta),
\end{equation}

except for $z$ in a set of discs, the sum of whose radii is at
most $\eta R$.\\

It remains to consider the term $\log \prod_{n =1}^m R/r_n$ in
(3.5). We have

\begin{eqnarray*}
 \sum_{n=1}^m \log  \frac{R}{r_n}
 & = & m\log R - \sum_{n=1}^m \log r_n\\
& = & m\log R - m \log r_m + \sum_{n=1}^{m-1} n(\log r_{n+1} -
\log
 r_n)\\
 & = & \int_{0}^{R} \frac{n(t)}{t}dt = N(R).
\end{eqnarray*}

 Together with (3.1), (3.4), (3.5) and (3.7) this gives the required result.
\end{proof}

\section{A new $\cos \pi \rho$-type theorem}
\setcounter{equation}{0}

Recall that Theorem~2 states that, for a transcendental entire
function $f$ of order less than $1$, there is a close
relationship between the minimum and maximum modulus of $f$ near
places where the maximum modulus of $f$ is small; in particular
if $f$ has order zero, this relationship exists on a set of upper
density 1. It is sufficient to prove Theorem~2 in the case that
$f(0) = 1$. Otherwise we have $g(0)=1$ for some function of
the form $g(z) = f(z)/(az^p)$, where $a \neq 0$, $p\in \N$, and the result for $f$ follows
from the result for $g$ by a routine calculation.\\

We prove Theorem 2, in the case that $f(0) = 1$, by using the bounds for the maximum and minimum modulus
given in Lemma~3.3 and Lemma~3.5. These lemmas show that we can
find values of $r$ for which the minimum and maximum modulus are
close by finding values of $r$ for which $Q(r)/N(r)$ is (relatively)
small. We do this by finding values of $r$ for which $a(r)/B(r)$
is small and then using Lemma~3.1 to deduce that, for such
values, $Q(r)/N(r)$ is also small.

\newtheorem{3.1}{Lemma}[section]
\begin{3.1}
Let $f$ be a transcendental entire function of order $\rho$, $\rho
<\alpha < 1$ with $f(0)=1$ and put
\[
  \epsilon(r) = \frac{\log \log M(r)}{\log r} \; \mbox{ and } \;
  \delta(r) = \max_{r \leq t \leq r^{1/(1-\alpha)}} \epsilon(t).
\]
Let $\mu,$ $\nu$ be functions such that, for large $r$,
\[
0< \mu(r) \leq 1,\; 0<\nu(r) \leq 1/4 \; \mbox{ and } \; \mu(r) \nu(r) >
\frac{2 \epsilon(r)}{1 - \delta(r)}.
\]
Then, for large $r$, there exists $R \in (r^{(1 - \delta(r))(1 -
\mu(r))},r^{1 - \delta(r)})$ such that
\[
 \frac{a(R)}{B(R)} \leq \nu(r)\; \mbox{ and hence }\; \frac{Q(R)}{N(R)} \leq 8\nu(r).
\]
\end{3.1}
\begin{proof} Suppose that there exist
arbitrarily large values of $r$ such that
\begin{equation}
\frac{a(t)}{B(t)} > \nu(r), \mbox{ for all } t \in (r^{(1 -
\delta(r))(1 - \mu(r))},r^{1 - \delta(r)}).
\end{equation}
Now $a(t) = tB'(t)$ and so it follows from (4.1) that
\begin{eqnarray*}
\int_{r^{(1 - \delta(r))(1 - \mu(r))}}^{r^{1 - \delta(r)}}
\frac{B'(t)}{B(t)}dt & > &
 \nu(r)\int_{r^{(1 - \delta(r))(1 - \mu(r))}}^{r^{1 -
 \delta(r)}}
 \frac{1}{t}dt\\
 & = & \nu(r) \log r^{(1 - \delta(r))\mu(r)}\\
 & = & \nu(r) \mu(r) (1 - \delta(r)) \log r\\
 & > & 2 \epsilon(r)\log r.
\end{eqnarray*}
Thus, for arbitrarily large values of $r$, we have
\begin{equation}
\log B(r^{1 - \delta(r)}) > 2 \epsilon(r)\log r.
\end{equation}

It follows from Lemma~3.2, however, that, for large $r$,
\begin{eqnarray*}
\log B(r^{1 - \delta(r)}) &  < &
 \log (3 \log M(r) \log r )
  <  \log ((\log M(r))^2)\\
  & = & 2 \log \log M(r) = 2 \epsilon(r) \log
  r.
\end{eqnarray*}

This contradicts (4.2) and so our original supposition must be
false. \\

Hence, for large $r$, there exists $R \in (r^{(1 - \delta(r))(1 - \mu(r))},r^{1 -
\delta(r)})$ such that
\[
 a(R) \leq \nu(r) B(R) \leq B(R)/4,
\]
and so, by Lemma~3.1(e),
\[
 B(R) \leq N(R) + 2a(R) \leq N(R) + B(R)/2.
\]
Thus, for such $R$,
\begin{equation}
B(R) \leq 2 N(R).
\end{equation}
We know from Lemma~3.1(d) that $Q(R) \leq 4a(R)$ and so, for such $R$,
\[
 \frac{Q(R)}{N(R)} \leq \frac{8a(R)}{B(R)} \leq 8 \nu(r),
\]
as required.
\end{proof}

{\it Proof of Theorem~2.}
 Let $f$ be a transcendental entire
function of order $\rho$, $\rho <\alpha < 1$ with $f(0)=1$ and put
\[
  \epsilon(r) = \frac{\log \log M(r)}{\log r} \; \mbox{ and } \;
  \delta(r) = \max_{r \leq t \leq r^{1/(1-\alpha)}} \epsilon(t).
\]
Let $0 < \eta < 1/2$ and let $\mu,$ $\nu$ be functions such
that, for large $r$,
\[
0<\mu(r) \leq 1,\; 0<\nu(r) \leq 1/4 \; \mbox{ and } \; \mu(r) \nu(r) >
\frac{2 \epsilon(r)}{1 - \delta(r)}\,.
\]

It follows from Lemma~3.3, Lemma~3.5 and Lemma~4.1 that, for large
$r$, there exists $R \in (r^{(1 - \delta(r))(1 - \mu(r))},r^{1 -
\delta(r)})$ such that
\begin{eqnarray*}
\frac{\log m(t)}{\log M(t)} & > & \frac{N(R) - (1 +
\log(2e/\eta))Q(R)}{N(R) + Q(R)}\\
& \geq & \frac{1 - (1 + \log(2e/\eta))8\nu(r)}{1 + 8
\nu(r)}\\
& > & 1 - 20\log(2e/\eta) \nu(r),
\end{eqnarray*}
for $0 \leq t \leq R/2$, except in a set of intervals, the sum of
whose lengths is at most $\eta R$. This completes the proof of
Theorem~2.

\section{A new growth condition for Baker's question}
\setcounter{equation}{0}

We begin this section by showing how Theorem~3 follows from
Theorems 1 and 2. Let $f$ be a transcendental entire function of
order zero and put
\[
  \epsilon(r) = \frac{\log \log M(r)}{\log r}\; \mbox{ and }\;
  \delta(r) = \max_{r \leq t \leq r^{2}} \epsilon(t).
\]

For large $r$, we have $0 < \epsilon(r) \leq 1/144$ and $\delta(r) < 1/3$, and so we can apply
Theorem~2 with $\alpha = 1/2$, $\eta = 1/8$, $\mu(r) =
\sqrt{\epsilon(r)}$ and $\nu(r) = 3\sqrt{\epsilon(r)}$.\\

 Now
let $R_1>0$ be such that the sequence $R_n$ defined by $R_{n+1} =
M(R_n)$ tends to $\infty$ and the sequence $\epsilon_n$ defined by
\[
  \epsilon_n = \max_{R_n \leq r \leq R_{n+1}} \epsilon(r)
\]
is positive, and suppose that $\sum_{n \in \N} \sqrt{\epsilon_n}
< \infty$. Let $r \in [R_n,R_{n+1})$, for some $n \in \N$. We may
assume that $R_1$ has been chosen sufficiently large to ensure
that, by Theorem~2, there exists
\[
t \in (r^{(1 - \delta(r))(1 - \sqrt{\epsilon(r)})}/4, r^{1 -
\delta(r)}) \subset (r^{(1 - \delta(r))(1 - \sqrt{\epsilon(r)}) -
2/\log r}, r^{1 - \delta(r)})
\] such that
\[
\log m(t)  >  \left(1 - 60 \log(16e) \sqrt{\epsilon(r)}\right)
\log M(t)
\]
\begin{equation}
 >  \left(1 - 230 \sqrt{\epsilon_n}\right) \log M(t).
\end{equation}

Now let $L>1$. We may assume that $R_{n+2}>R_{n+1}^{2L}$ and so
$r^{2L} \in [R_n,R_{n+2})$. Thus
\[
\epsilon(r^L) \leq \delta(r^L) \leq \delta_n, \mbox{ where } \delta_n = \max\{\epsilon_n,
\epsilon_{n+1}\}.
\]
We may also assume that $R_1$ has been chosen sufficiently large to ensure
that $\epsilon_n < 1$ for all $n \in \N$ and that, by (5.1) and
Lemma~2.2, there exists
\[
t \in (r^{L(1 - \delta_n)(1 - \sqrt{\delta_n}) - 2/\log R_n},
r^{L}) \subset (r, r^L)
\]
 such that
\begin{eqnarray*}
\log m(t) & > & \left( 1 - 230 \sqrt{\delta_n}\right) \log M\left(
r^{L(1 - \delta_n)(1 -
\sqrt{\delta_n}) - 2/\log R_n}\right)\\
& > & L\left(1 - 230 \sqrt{\delta_n}\right)\left((1 - \delta_n)(1
- \sqrt{\delta_n}) -
\frac{2}{\log R_n}\right) \log M(r)\\
& > & L\left(1 - 232 \sqrt{\delta_n} - \frac{2}{\log R_n}\right)
\log M(r).
\end{eqnarray*}

Since $\log R_n
> 2^n$, for large $n \in \N$, we have $\sum_{n \in \N}
\frac{1}{\log R_n} < \infty$. Thus
\[
 \sum_{n \in \N} \left( 232 \sqrt{\delta_n} + \frac{2}{\log R_n} \right)
  < \sum_{n \in \N} \left( 464 \sqrt{\epsilon_n} + \frac{2}{\log
 R_n}\right) <
 \infty
\]
and so it follows from Theorem~1 with $a_n = 232 \sqrt{\delta_n}
+ 2/ \log R_n$ that $F(f)$ has no unbounded
components. This completes the proof of Theorem~3.\\

We now use the result of Theorem~3 to prove Theorem~4.\\

{\it Proof of Theorem~4.} Suppose that there exist $m \in \N$ and
$R>0$ such that
\[
 \log \log M(r) < \frac{\log r}{\log^m r}, \mbox{ for }r>R.
\]
Now take $R_1>R$ so large that the sequence $R_n$ defined by
$R_{n+1} = M(R_n)$ tends to $\infty$ and the sequence
$\epsilon_n$ defined by
\[
  \epsilon_n = \max_{R_n \leq r \leq R_{n+1}} \frac{\log \log M(r)}{\log
  r}
\]
is positive. Then
\begin{equation}
 \epsilon_n \leq \frac{1}{\log^m R_n}.
\end{equation}
We will show that $\sum_{n \in \N} \sqrt{\epsilon_n} < \infty$
and hence, by Theorem~3, $F(f)$ has no unbounded components. To
do this, we show that the values of $n$ satisfying
\begin{equation}
 \epsilon_n \geq 1/n^4
\end{equation}
are relatively sparsely distributed.\\

We begin by noting that if $n$ satisfies (5.3) then there exists
$r_n \in [R_n, R_{n+1}]$ with $M(r_n) \geq \exp(r_n^{1/n^4})$ and
so
\begin{equation}
R_{n+2} = M(R_{n+1}) \geq M(r_n) \geq \exp(r_n^{1/n^4}) \geq
\exp(R_n^{1/n^4}).
\end{equation}
Suppose that $N$ satisfies (5.3) and that there exist $N_1,N_2$
satisfying (5.3) with
\begin{equation}
N_1 \geq N+2 \mbox{ and } e^N \geq N_2 \geq N_1 + 2.
\end{equation}
Then, by (5.4) and (5.5),
\[
 R_{N_2} \geq M(R_{N_1 + 1}) \geq \exp(R_{N_1}^{1/N_1^4}) \geq \exp(R_{N_1}^{1/e^{4N}})
\]
and
\[
 R_{N_1} \geq M(R_{N+1}) \geq \exp(R_N^{1/N^4}).
\]

Thus
\[
  R_{N_2} \geq \exp(\exp(R_{N}^{1/N^4}/e^{4N})).
\]
We claim that if, in addition, $N$ is sufficiently large, then
\begin{equation}
R_{N_2} \geq e^{R_N}.
\end{equation}
This is true since
\begin{eqnarray*}
\exp(\exp(R_{N}^{1/N^4}/e^{4N})) & \geq & e^{R_N}\\
\iff \frac{\log R_N}{N^4} - 4N & \geq & \log \log R_N\\
\iff \frac{\log R_N}{\log \log R_N}  & \geq & N^4\left(
\frac{4N}{\log \log R_N}+ 1 \right).
\end{eqnarray*}
Now, for large $N$, we have $\log R_N \geq 2^N$. Also, $r/\log r$
is increasing for $r \geq e$, so
\[
  \frac{\log R_N}{\log \log R_N} \geq \frac{2^N}{N \log 2} \geq
  N^4\left(
\frac{4N}{\log \log R_N}+ 1 \right),
\]
for large $N$. Thus (5.6) is true for large $N$, $N_1$ and $N_2$ satisfying (5.3) and (5.5).\\

Now suppose that $N$ satisfies (5.3) and there exist $N_i$, $1 \leq i \leq 2m$,
satisfying (5.3), with $e^N \geq N_{i+1} \geq N_i + 2$, for $1
\leq i \leq 2m-1$, and $N_1 \geq N + 2$. It follows from (5.2) and repeated application of
(5.6) that, if $N$ is sufficiently large and $N_{2m} \leq n \leq
e^N$, then
\begin{equation}
\epsilon_n \leq \frac{1}{\log^m R_n} \leq \frac{1}{\log ^m
R_{N_{2m}}} \leq \frac{1}{\log^m \exp^m R_N} = \frac{1}{R_N} \leq
\frac{1}{\exp(2^N)} \leq \frac{1}{n^4}.
\end{equation}
So, for large $N$, there are at most $4m+1$ values of $n$ with $N
\leq n \leq e^N$ such that $n$ satisfies (5.3). Thus, by (5.2),
\[
\sum_{n=N}^{[e^N]} \sqrt{\epsilon_n} < \frac{4m+1}{\sqrt{\log ^m
R_N}} + \sum_{n=N}^{[e^N]} \frac{1}{\sqrt{n^4}}.
\]
Now, for large $N$, we have $R_N > N$ and so there exist constants
$N_0 \in \N$ and $C>0$ such that
\[
\sum_{n\in \N} \sqrt{\epsilon_n} < C + \sum_{n \in \N}
\frac{1}{n^2} + (4m+1) \sum_{k \in \N} \frac{1}{\sqrt{\log^m (\exp^k N_0)}}
< \infty.
\]
This completes the proof of Theorem~4.

\section{Examples}
\setcounter{equation}{0}

 Let
\[
  f(z) = \prod_{m=1}^{\infty} \left(1 -
  \frac{z}{r_m}\right)^{r_m^{\epsilon_m}},
\]
where $\epsilon_m$ is a decreasing null sequence of positive
terms with $\epsilon_1 < 1$ and $r_m$ is a strictly increasing
sequence with $r_1 \geq 4$ such that $r_m^{\epsilon_m} \in \N$
for $m \in \N$. We will show that if $\epsilon_m$ and $r_m$
satisfy certain conditions, then $f$ is a function of order zero
that fails to satisfy the hypotheses of Theorem~1. In particular,
we give conditions which are sufficient to ensure that $f$ is a
function of order zero that fails to satisfy Hinkkanen's
condition (1.3). We begin by giving a condition on $r_k$ and
$\epsilon_k$ which is sufficient to ensure that $f$ has order
zero.

\newtheorem{6.1}{Lemma}[section]
\begin{6.1}
Suppose that
\begin{equation}
 r_{k+1}^{\epsilon_{k+1}} \geq r_k^2, \mbox{ for all } k \in \N.
\end{equation}
Then $f$ has order zero.
\end{6.1}

\begin{proof} Let $k \in \N$ and let $r_k \leq r < r_{k+1}$. Then,
by (6.1),
\begin{equation}
 M(r) = \prod_{m=1}^{\infty} \left( 1 + \frac{r}{r_m}\right)^{r_m^{\epsilon_m}}
 \leq r^{\sum_{m=1}^kr_m^{\epsilon_m}} \left( 1 + \frac{r}{r_{k+1}}\right)^{r_{k+1}^{\epsilon_{k+1}}}
 \prod_{m=k+2}^{\infty} \left( 1 + \frac{1}{r_m^{1/2}}\right)^{r_m^{\epsilon_m}}.
\end{equation}
Now, for $1 \leq m \leq k-1$, it follows from (6.1) that
\[
 r_m^{\epsilon_m} < r_m \leq r_{k-1} \leq r_k^{\epsilon_k/2}.
\]
Also, by (6.1), for $k\geq 2$ we have $r_k^{\epsilon_k} \geq 4^k
\geq (k-1)^2$ and so
\[
  \sum_{m=1}^kr_m^{\epsilon_m} \leq r_k^{\epsilon_k} +
  (k-1)r_k^{\epsilon_k/2}\leq 2r_k^{\epsilon_k}.
\]

So, for large $k$, it follows from (6.1) and (6.2) that
\begin{eqnarray*}
\log M(r) & \leq & \sum_{m=1}^kr_m^{\epsilon_m} \log r +  \frac{
r_{k+1}^{\epsilon_{k+1}}r}{r_{k+1}} +
\sum_{m=k+2}^{\infty}r_m^{\epsilon_m - 1/2} \\
& \leq & 2r_k^{\epsilon_k} \log r + r^{\epsilon_{k+1}} + \sum_{m=k+2}^{\infty}r_m^{-1/4}\\
& \leq & 2r_k^{\epsilon_k} \log r + r^{\epsilon_{k+1}} + 1,
\end{eqnarray*}
since $r_m \geq 4^m$. Thus, for large $k$,
\begin{equation}
\log M(r) \leq 3  r^{\epsilon_k} \log r, \mbox{ for } r_k \leq r
< r_{k+1}.
\end{equation}
So, for large $k$,
\[
 \frac{\log \log M(r)}{\log r} \leq \frac{\log(3 r^{\epsilon_k}\log r)}{\log
 r} = \epsilon_k + \frac{\log(3 \log r)}{\log r}, \mbox{ for } r_k \leq r <
 r_{k+1},
\]
and hence $\lim_{r \to \infty} \frac{\log \log M(r)}{\log r} = 0$.
\end{proof}

The next lemma is the main result in this section.

\newtheorem{6.2}[6.1]{Lemma}
\begin{6.2}
Suppose that
\begin{equation}
 r_{k+1}^{\epsilon_{k+1}} \geq  r_k^{k+1}, \mbox{ for all } k \in \N.
\end{equation}
Then, given any $L > 1$, there exists $K_L$ such that, for all
$k \geq K_L$,
\[
 \frac{\log m(t)}{\log M(r_k^{1/L})} \leq L(1 - \epsilon_k/4), \mbox{ for
 all } t \in (r_k^{1/L},r_k).
\]
\end{6.2}
\begin{proof} We fix $L > 1$. First note that it follows from
(6.4) that, for large $k$ (depending on $L$),
\[
M(r_k^{1/L})  =  \prod_{m=1}^{\infty}\left( 1 +
\frac{r_k^{1/L}}{r_m}\right) ^{r_m^{\epsilon_m}}
> \prod_{m=1}^{k-1}\left( \frac{r_k^{1/L}}{r_m}\right) ^{r_m^{\epsilon_m}}
\geq \prod_{m=1}^{k-1} r_k^{(1/L)(1 -
\epsilon_k/8)r_m^{\epsilon_m}}.
\]
Thus, for large $k$,
\begin{equation}
 M(r_k^{1/L}) \geq r_k^{(1/L)(1 - \epsilon_k/8)N_k}, \mbox{ where }  N_k
 = \sum_{m=1}^{k-1}r_m^{\epsilon_m}.
\end{equation}

Now, take $t \in (r_k^{1/L},r_k)$ for some $k \in \N$. If $k$ is large enough, then $r_{k-1} < r_k^{1/L}$, so
\begin{eqnarray*}
m(t) & = & \prod_{m=1}^{k-1} \left( \frac{t}{r_m}- 1
\right)^{r_m^{\epsilon_m}} \prod_{m=k}^{\infty} \left( 1 -
\frac{t}{r_m} \right)^{r_m^{\epsilon_m}}\\
& \leq & \prod_{m=1}^{k-1} t^{r_m^{\epsilon_m}} \left( 1 -
\frac{t}{r_k} \right)^{r_k^{\epsilon_k}}.
\end{eqnarray*}
Thus, for large $k$,
\begin{equation}
 m(t) \leq t^{N_k} \left( 1 -
\frac{t}{r_k} \right)^{r_k^{\epsilon_k}}, \mbox{ for }t \in
(r_k^{1/L},r_k).
\end{equation}

Differentiating, we see that the right-hand side of this
inequality takes its maximum value when
\[
 \frac{N_k}{t} - \frac{r_k^{\epsilon_k}}{r_k(1 - t/r_k)} = 0;
\]
that is, at $t_0$ where
\[
 t_0 = N_kr_k^{1 - \epsilon_k}(1 - t_0/r_k).
\]
Thus, by (6.6), for large $k$ and for all $t \in (r_k^{1/L},r_k)$,
\[
 m(t) \leq t_0^{N_k}
 \leq N_k^{N_k}r_k^{(1-\epsilon_k)N_k},
\]
and so, by (6.5),
\begin{equation}
\frac{\log m(t)}{\log M(r_k^{1/L})} \leq \frac{N_k \log N_k + N_k(1 -
\epsilon_k) \log r_k}{(N_k/L)(1 - \epsilon_k/8) \log r_k} =
L \, \frac{\log N_k + (1 - \epsilon_k) \log r_k}{(1 - \epsilon_k/8)
\log r_k}.
\end{equation}
It follows from (6.4) that
\[
N_k = \sum_{m=1}^{k-1}r_m^{\epsilon_m} \leq (k-1)
r_{k-1}^{\epsilon_{k-1}} \leq r_{k-1}^2 < r_k^{2\epsilon_k/k}
\]
and so, for large $k$,
\[
\log N_k + (1 - \epsilon_k) \log r_k < (2\epsilon_k/k) \log r_k + (1
- \epsilon_k) \log r_k < (1 - \epsilon_k/2) \log r_k.
\]
Thus, by (6.7), for large $k$ and for all $t \in (r_k^{1/L},r_k)$,
\[
 \frac{\log m(t)}{\log M(r_k^{1/L})} \leq L\frac{(1 - \epsilon_k/2) \log r_k}{(1 - \epsilon_k/8) \log r_k}
 < L(1 - \epsilon_k/4),
\]
as required.
\end{proof}

The next lemma gives a condition on $r_k$ and $\epsilon_k$
which is sufficient to ensure that if $L > 1$ then, for large values of $n$, each interval of the form $[R_n, R_{n+1})$ contains at most one point of the form $r_k^{1/L}$.

\newtheorem{6.3}[6.1]{Lemma}
\begin{6.3}
Suppose that
\begin{equation}
 r_{k+1}^{\epsilon_{k+1}} \geq e^{r_k}, \mbox{ for all } k \in \N,
\end{equation}
and let $R_1>0$ be such that the sequence $R_n$ defined by $R_{n+1} = M(R_n)$ tends to $\infty$. Now fix $L > 1$ and, for large $k$, let $R_{n_k}$ be
such that
\[
 r_k^{1/L}
\in [R_{n_k}, R_{n_k + 1}), \mbox{ where } R_{n+1} = M(R_n).
\]
Then, for large $k$, $R_{n_{k+1}} \geq R_{n_k + 1}$ and hence $n_{k+1} > n_k$.
\end{6.3}
\begin{proof}
 Note that
 \[
 R_{n_k + 1} = M(R_{n_k}) \leq M(r_k^{1/L}) < M(r_k).
 \]
 Since (6.8) is satisfied, (6.1) and hence (6.3) are also satisfied. It follows from (6.3) and (6.8) that,
  for large $k$,
  \[
    R_{n_k+1} < M(r_k) \leq \exp(3 r_k^{\epsilon_k} \log r_k) < \exp(r_k) <
    r_{k+1}^{1/L}.
  \]
  Thus $R_{n_{k+1}} \geq R_{n_k + 1}$.
\end{proof}

To construct the required example, we consider the function $f$
given by (1.7), where $\epsilon_1<1$ and
\begin{equation}
\sum_{m=1}^{\infty} \epsilon_m = \infty,
\end{equation}
and the sequence $r_m$ is chosen to satisfy $r_1\ge 4$ and (6.8),
and hence satisfies (6.1) and (6.4) (because $e^x \geq x^{k+1}$,
for $x \geq (k+1)^2$ and $r_k \geq 4^k \geq (k+1)^2$). In
particular, $f$ has
order zero by Lemma~6.1.\\

Now suppose that $a_n$ is a positive sequence for which condition
(1.4) of Theorem~1 is satisfied and let $R_1>0$ be such that the
sequence $R_n$ defined by $R_{n+1} = M(R_n)$ tends to $\infty$.
Then for $k$ large enough, there exists $t_k \in (r_k^{1/L},r_k)$
such that
\[
\frac{\log m(t_k)}{\log M(r_k^{1/L})} \geq L(1 - a_{n_k}),
\]
where $n_k$ is the integer such that
\[
 r_k^{1/L}
\in [R_{n_k}, R_{n_k + 1}).
\]
Thus, by Lemma~6.2, we have
\[
 L(1 - \epsilon_k/4) \geq L(1 - a_{n_k}), \mbox{ so } a_{n_k} \geq \epsilon_k/4,
\]
for all large $k$. Therefore, by Lemma~6.3 and (6.9), we deduce that $ \sum_{n
\in \N} a_n = \infty$, so $f$ is a transcendental entire function
of order zero for which there is no positive sequence $a_n$
satisfying all
the conditions of Theorem~1.\\

{\it Remark.} Note that this function cannot satisfy the growth restriction
in Theorem~4 since this would imply that it satisfies all the conditions of Theorem 1.\\

Finally, suppose that
\begin{equation}
\epsilon_k (\log r_k)^{1/k} \to \infty,\; \mbox{ as } k \to \infty.
\end{equation}
Then, given $L>1$ and $C,\delta > 0$, we have
\[
  \epsilon_k/4 > \frac{C}{(\log r_k^{1/L})^{\delta}}, \mbox{ for
  large } k,
\]

and so, if (6.4) is also satisfied, then it follows from Lemma~6.2
that, for large $k$,

\[
\frac{\log m(t)}{\log M(r_k^{1/L})} \leq L\left( 1 -
\frac{C}{(\log r_k^{1/L})^{\delta}}\right), \mbox{ for
 all } t \in (r_k^{1/L},r_k).
\]

If (6.4) is satisfied, then so is (6.1) and so, by Lemma~6.1, $f$
has order zero. Thus, if (6.10) and (6.4) are satisfied, then $f$
is a function of order zero that fails to satisfy Hinkkanen's
condition (1.3).\\

\section{Regularity conditions}
\setcounter{equation}{0}

Many authors have shown that a function of order $\rho < 1/2$ has
no unbounded Fatou components if the growth is sufficiently
regular. We now prove Theorem~5 and show that the regularity
condition
(1.8) includes various known regularity conditions.\\

We use the following result which is due to Baker [\ref{B58},
Satz 1]. (This result also follows from the version of the $\cos
\pi \rho$ theorem proved by Barry [\ref{Ba}].)

\newtheorem{4.1}{Lemma}[section]
\begin{4.1}
Let $f$ be a transcendental entire function of order $\rho < 1/2$.
There exist $m>1$ and $R>0$ such that, for all $r>R$, there
exists $r'$ satisfying
\[
 r \leq r' \leq r^m, \; \mbox{ with } m(r') = M(r).
\]
\end{4.1}

The following result is a direct consequence of Lemma~7.1.

\newtheorem{4.4}[4.1]{Lemma}
\begin{4.4}
Let $f$ be a transcendental entire function of order $\rho < 1/2$
and let $m$ be the constant given in Lemma~7.1. Suppose that there
exist a sequence $R_n\to\infty$ defined by $R_{n+1} = M(R_n)$ and
a function $\psi: [R_1, \infty) \to [R_1,\infty)$ such that
$\psi(r) \geq r$ and
\begin{equation}
 M(\psi(R_n)) \geq \psi(R_{n+1})^m, \; \mbox{ for } n \in \N.
\end{equation}
Then there exists a sequence $\rho_n$ such that
\begin{enumerate}
\item $R_n \leq \psi(R_n) \leq \rho_n \leq \psi(R_n)^m$,

\item $m(\rho_n) \geq \psi(R_{n+1})^m$.
\end{enumerate}
\end{4.4}
\vspace{10pt}
Theorem~5 follows easily from Lemma~7.2.\\

{\it Proof of Theorem~5.}\;
If the hypotheses of Theorem~5 are satisfied, then so are the
hypotheses of Lemma~7.2. This implies that the
 hypotheses of Lemma~2.1 are also satisfied, since $\psi(R_n)^m =
R_n^{c(n)}$, where $c(n)>1$. Thus $F(f)$ has no unbounded
components. This completes the proof of Theorem~5.\\

We stated earlier that many of the known regularity conditions
associated with Baker's question can be written in the form (1.8).
For example, we showed in [\ref{S}, Theorem~C] that a
transcendental entire function of order $\rho < 1/2$ has no
unbounded Fatou components if there exists a finite constant $c$
such that
\begin{equation}
\frac{\log M(2r)}{\log M(r)} \to c\; \mbox{ as } r \to \infty.
\end{equation}
We showed in [\ref{S}, equation (4.2)] that, for such a function
with $\rho > 0$, the inequality (1.8) is satisfied with $\psi(r) = r^2$;
for functions of order zero that satisfy (7.2) a more delicate
argument
was required.\\

Another regularity condition was given by Anderson and Hinkkanen
[\ref{AH}]. They showed that an entire function of order $\rho <
1/2$ has no unbounded Fatou components if there exists $c>0$ such
that the function $\phi(x) = \log M(e^x)$ satisfies
\begin{equation}
  \frac{\phi'(x)}{\phi(x)} \geq \frac{1+c}{x}, \; \mbox{ for large } x.
\end{equation}
Note that (7.3) implies that, for large $x$ and all $k > 1$,
\[
 \int_x^{kx} \frac{\phi'(t)}{\phi(t)}\,dt  \geq \int_x^{kx} \frac{1+c}{t}\,dt
\]
and so
\[
 \log (\phi(kx)/\phi(x)) \geq (1+c)\log k;
\]
that is,
\[
  \phi(kx) \geq k^{1+c} \phi(x).
\]
Thus, if $m>1$ is given, then by taking $k = m^{1/c} > 1$, we
obtain
\[
 \phi(kx) \geq km \phi(x), \mbox{ for large } x,
\]
and so
\[
 M(r^k) \geq  M(r)^{km}, \; \mbox{ for large } r.
\]
Hence, if (7.3) holds, then the conditions of Theorem~5 are
satisfied with
$\psi(r) = r^k$.\\

{\it Remark.} We mention here the paper of Hua and Yang [\ref{HY}] in which further such
regularity conditions are stated. Unfortunately, the proofs of the main results in [\ref{HY}]
appear to contain gaps, as pointed out in the survey article [\ref{H}].\\

Finally, we prove Theorem~6.\\

{\it Proof of Theorem~6.}\;
 Let $f$ be a transcendental entire
function of order $\rho < 1/2$ and suppose that there exist $n
\in \N$ and $0<q<1$  such that
\[
 M(r) \geq \exp^{n+1}((\log^nr)^q), \; \mbox{ for large } r.
\]
 Now let $m>1$ and let
\[
 \psi(r) = \exp^{n}((\log r)^p), \; \mbox{ where } pq>1.
\]
We will show that, with $M(r)$ and $\psi(r)$ as given above,
for large $r$ we have $M(\psi(r)) \geq (\psi(M(r)))^m$. Since
$\psi(r) \geq r$ for $r \geq e$, this is sufficient to show that
the conditions of Theorem~5 are satisfied. The result then
follows.\\

We begin by noting that, for large $r$,
\[
M(\psi(r)) \geq \exp^{n+1}((\log^n(\exp^n((\log r)^p)))^q) =
\exp^{n+1}((\log r)^{pq}).
\]
Also, since $f$ has order less than 1/2, we have, for large $r$,
\[
 \psi(M(r))^m \leq \psi(e^{r^{1/2}})^m =
 \exp^n(r^{p/2})^m.
\]
So it is sufficient to show that, for large $r$,
\[
 \exp^{n+1}((\log r)^{pq}) \geq \exp^n(r^{p/2})^m;
\]
that is,
\begin{equation}
\exp^{n-1}((\log r)^{pq}) \geq \log m + \exp^{n-2}(r^{p/2}), \;
\mbox{ if } n \geq 2,
\end{equation}
or
\begin{equation}
(\log r)^{pq} \geq \log m + (p/2) \log r, \; \mbox{ if } n = 1.
\end{equation}

Since (7.5) is clearly true for large $r$, it remains to show that
(7.4) is also true for large $r$. We note that, for large $r$,
\begin{eqnarray*}
\exp^{n-1}((\log r)^{pq}) & = & \exp^{n-2}(\exp((\log r)^{pq}))\\
& \geq & \exp^{n-2}(2 r^{p/2})\\
& \geq & 2 \exp^{n-2}(r^{p/2})\\
& \geq & \log m + \exp^{n-2}(r^{p/2})
\end{eqnarray*}
and so (7.4) is true, as required. This completes the proof of
Theorem~6.

\section*{References}
\begin{enumerate}

\item \label{AH} J.M. Anderson and A. Hinkkanen. Unbounded domains
of normality. {\it Proc. Amer. Math. Soc. 126 (1998), 3243--3252.}

\item \label{B58} I.N. Baker. Zusammensetzungen ganzer Funktionen.
{\it Math. Z. 69 (1958), 121--163.}

\item \label{B81} I.N. Baker. The iteration of polynomials and
transcendental entire functions. {\it J. Austral. Math. Soc.
(Series A) 30 (1981), 483--495.}

\item \label{Ba} P.D. Barry. On a theorem of Besicovitch.
{\it Quart. J. Math. Oxford Ser. (2) 14 (1963), 293--302}.

\item \label{Berg} W. Bergweiler. Iteration of meromorphic functions.
{\it Bull. Amer. Math. Soc. 29 (1993), 151--188.}

\item \label{BH} W. Bergweiler and A. Hinkkanen. On
semiconjugation of entire functions. {\it Math. Proc. Camb. Phil.
Soc. 126 (1999), 565--574.}

\item \label{Ca} H. Cartan. Sur les syst\`{e}mes de fonctions holomorphes
\`{a} vari\'{e}t\'{e}s lin\'{e}aires et leurs applications. {\it
Ann. Sci. \'{E}cole Norm. Sup. (3) 45 (1928), 255--346}.

\item \label{C} M.L. Cartwright. {\it Integral functions.}
Cambridge Tracts in Mathematics and Mathematical Physics, No. 44,
Cambridge University Press, 1962.

\item \label{E} A.E. Eremenko. On the iteration of entire
functions. {\it Dynamical systems and ergodic theory,} Banach
Center Publ. 23 (Polish Scientific Publishers, Warsaw, 1989)
339--345.

\item \label{H2} A. Hinkkanen. Entire functions with no unbounded
Fatou components. {\it Complex analysis and dynamical systems II,
217--226.} Contemp. Math. 382, Amer. Math. Soc., Providence, RI,
2005.

\item \label{H} A. Hinkkanen. Entire functions with bounded Fatou
components. To appear in {\it Transcendental dynamics and complex
analysis.}  Cambridge University Press, 2008.

\item \label{HM} A. Hinkkanen and J. Miles. Growth conditions for entire functions with only bounded Fatou components. {\it Preprint.}

\item \label{HY} X. Hua  and C.C. Yang. Fatou components of entire
functions of small growth {\it Ergodic Theory Dynam.\ Systems 19
(1999),  1281--1293.}

\item \label{Zip} P.J. Rippon and G.M. Stallard. On sets where
iterates of a meromorphic function zip towards infinity. {\it
Bull. London Math. Soc. 32 (2000), 528--536.}

\item \label{RS05} P.J. Rippon and G.M. Stallard. On questions of
Fatou and Eremenko. {\it Proc. Amer. Math. Soc. 133 (2005),
1119--1126.}

\item \label{RS} P.J. Rippon and G.M. Stallard. Escaping points of entire function of small growth. {\it Preprint.}

\item \label{Si} A.P. Singh.  Composite entire functions with no unbounded
Fatou components.  {\it J. Math. Anal. Appl., 335 (2007),
907--914.}

\item \label{S} G.M. Stallard. The iteration of entire functions
of small growth. {\it Math. Proc. Camb. Phil. Soc. 114 (1993),
43--55.}

\item \label{T} E.C. Titchmarsh. {\it The theory of functions.}
Oxford University Press, 1939.

\item \label{W} Y. Wang. Bounded domains of the Fatou set of an
entire function. {\it Israel J. Math. 121 (2001), 55--60.}

\item \label{Z} Jian-Hua Zheng. Unbounded domains of normality of
entire functions of small growth. {\it Math. Proc. Camb. Phil.
Soc. 128 (2000), 355--361.}

\end{enumerate}

Department of Mathematics,

The Open University,

Walton Hall,

Milton Keynes,

MK7 6AA.

E-mail: p.j.rippon@open.ac.uk; g.m.stallard@open.ac.uk

\end{document}